\documentclass[10pt]{amsart}
\usepackage{amssymb,amscd}

\newcommand{\FF}{\mathbb F}

\newcommand{\CC}{\mathcal C}
\newcommand{\OO}{\mathcal O}

\newcommand{\NN}{\mathbb N}

\newcommand{\Dtil}{\widetilde{D}}
\newcommand{\Etil}{\widetilde{E}}

\newcommand{\set}[1]{\left\{#1\right\}}
\newcommand{\To}{\longrightarrow}

\newcommand{\isom}{\cong}

\DeclareMathOperator{\Aut}{Aut}

\DeclareMathOperator{\Hom}{Hom} 
 
\DeclareMathOperator{\Irr}{Irr}
\DeclareMathOperator{\coker}{coker}
\DeclareMathOperator{\core}{core} \DeclareMathOperator{\Res}{Res}

\DeclareMathOperator{\im}{im}
\DeclareMathOperator{\res}{res}

\DeclareMathOperator{\MOD}{mod}

\swapnumbers   
\theoremstyle{plain}

\newtheorem{thm}{Theorem}[section] 

\newtheorem{cor}[thm]{Corollary}

\newtheorem{prop}[thm]{Proposition}

\theoremstyle{definition}
\newtheorem{defn}[thm]{Definition}

\newtheorem{rem}[thm]{Remark}

\theoremstyle{remark}

\numberwithin{equation}{section} 


\renewcommand{\setminus}{\smallsetminus}
\renewcommand{\emptyset}{\varnothing}
\newcommand{\subsneq}{\varsubsetneqq}

\newbox\mybox

\def\arrover#1{\mathrel{
\setbox\mybox=\hbox spread 1.4em{\hfil$\scriptstyle#1$\hfil}
\vbox{\offinterlineskip\copy\mybox \hbox
to\wd\mybox{\rightarrowfill}}}}

\def\larrover#1{\mathrel{
\setbox\mybox=\hbox spread 1.4em{\hfil$\scriptstyle#1$\hfil}
\vbox{\offinterlineskip\copy\mybox \hbox
to\wd\mybox{\leftarrowfill}}}}

\def\ontoover#1{\mathrel{
\setbox\mybox=\hbox spread 1.4em{\hfil$\scriptstyle#1$\hfil}
\vbox{\offinterlineskip\copy\mybox \hbox
to\wd\mybox{\rightarrowfill\hskip-2.8mm $\rightarrow$}}}}

\def\leftontoover#1{\mathrel{
\setbox\mybox=\hbox spread 1.4em{\hfil$\scriptstyle#1$\hfil}
\vbox{\offinterlineskip\copy\mybox \hbox
to\wd\mybox{$\leftarrow$\hskip-2.8mm \leftarrowfill}}}}

\begin{document}

\title[On $p$-rank representations]{On $p$-rank representations}
\author{Nicolas Stalder}
\date{\today}
\address{Department of Mathematics, ETH Z\"{u}rich, CH-8092 Z{\"u}rich, Switzerland}
\subjclass[2000]{14H37; 20C20, 14F20}
\keywords{$p$-rank, Galois module}
\email{nicolas.stalder@math.ethz.ch}


\begin{abstract}
The $p$-rank of an algebraic curve $X$ over an algebraically closed field $k$
of characteristic $p>0$ is the dimension of the vector space $H^1_{et}(X,\FF_p)$.
We study the representations of finite subgroups $G \subset\Aut(X)$
induced on $H^1_{et}(X,\FF_p)\otimes k$, and obtain two main results:

First, the sum of the \emph{nonprojective} direct summands of the representation,
i.e. its \emph{core}, is determined explicitly by local data given
by the fixed point structure of the group acting on the curve. As a corollary,
we derive a congruence formula for the $p$-rank.

Secondly, the multiplicities of the \emph{projective} direct summands of quotient curves,
i.e. their \emph{Borne invariants}, are calculated in terms
of the Borne invariants of the original curve and ramification data. In
particular, this is a generalization of both Nakajima's equivariant
Deuring-Shafarevich formula and a previous result
of Borne in the case of free actions.
\end{abstract}


\maketitle
\tableofcontents


\section{Introduction}

We fix an irreducible, smooth and complete curve $X$ over an algebraically
closed field $k$ of positive characteristic $p$. The etale cohomology
group $H^1(X_{et},\FF_p)$ is a finite-dimensional vector space over
$\FF_p$. Its dimension, the \emph{$p$-rank} of $X$, is a global invariant of the curve.

If we fix a finite group $G$ of automorphisms of the curve $X$, then
$H^1(X_{et},\FF_p)$ becomes a finite-dimensional representation of $G$ over
$\FF_p$. Moreover, $H^1(X_{et},\FF_p)\otimes_{\FF_p}k$ is a finite-dimensional
representation of $G$ over $k$.

First results on determining this representation up to isomorphism
by local invariants of the curve and of the group action have been obtained
by Shoichi Nakajima \cite{Nak2}, under the assumption that $G$ is a $p$-group,
and by Niels Borne \cite{Bor}, under the assumption that $G$ operates
without fixed points. We continue this tradition, with no assumptions on
either the group $G$ or its action.

The local invariants cannot determine the representation completely, as
the example of an elliptic curve $E$ over
a field of characteristic $\neq 2$ shows: Such a curve always allows an
automorphism of order $2$, which stabilizes exactly $4$ points (and the
projection to the quotient curve is \emph{tamely} ramified in these points).
Namely, if the curve is given by the equation $y^2=f(x)$, consider the mapping
given by $(x,y)\mapsto (x,-y)$. However,
the $p$-rank of $E$ can be $0$ \emph{or} $1$, depending on whether this
curve is supersingular or not.

Accordingly, our results must be incomplete. Using the language of
modular representation theory, what we do determine completely is
the core of the representation (i.e., its ``non-semisimple'' part,
\emph{cf.} Section 1); this is the content of
Theorem~\ref{theoremA}. In a sense, this result is surprising,
since generally the non-semisimplicity of representations is what
makes modular representation theory more difficult than
representation theory in characteristic zero. Now the
representation is determined completely by its core and the
multiplicities of the indecomposable projective summands, which we
call \emph{Borne invariants} of the curve and introduce in Section
$5$. However, it is impossible to determine these by \emph{local}
invariants, as the above example shows.

The content of Theorem \ref{theoremB} is to determine explicitly in terms
of local data the Borne invariants of quotient curves $X/N$ with respect
to the quotient group $G/N$, for any normal subgroup $N\subset G$. This
gives a procedure for calculating the Borne invariants for those representations
of $G$ induced by quotient groups, in terms of the ``smaller'' curve $X/N$, and may
thus be regarded as a partial solution to the problem of determining
Borne invariants. In particular, if $N$ is a $p$-group this approach
gives all Borne invariants, and if $G$ itself is a $p$-group we recover
Nakajima's equivariant Deuring-Shafarevich formula.

This article is a digest of my diploma thesis at ETH Z\"urich; it is a
pleasure to thank my advisor Richard Pink for his support and guidance.
In particular, the proof of Proposition \ref{calcinv} is entirely due to him.



\section{Modular representation theory of finite groups}

It is customary to call a (finite dimensional) representation of a (finite)
group a \emph{modular representation} if the characteristic of
the field divides the order of the group. In this situation, the notions
of \emph{simple} and \emph{indecomposable} module no longer coincide,
as would be the case in characteristic $0$ by Maschke's theorem. This makes for
a richer representation theory, which we will now review.

In the following we shall fix an algebraically closed field
$k$ of characteristic $p>0$, a finite group $G$, and denote by $k[G]$ the group ring
of $G$ over $k$. All modules under consideration will be finitely generated left $k[G]$-modules,
and we identify finite dimensional representations of $G$ over $k$ with such modules.
All homomorphisms are assumed to be $k[G]$-linear.

\begin{defn}
A representation is \emph{simple} (or \emph{irreducible}) if it is nontrivial
and has no proper submodules.
We denote the set of isomorphism classes of simple modules by $\Irr G$. A representation  is
\emph{indecomposable} if it is nontrivial and admits no proper direct summands. It is \emph{projective}
if the functor $\Hom(P,-)$ is exact.
\end{defn}

\begin{thm}[Krull-Schmidt]
If $M$ is a representation, and $M\isom\oplus_{i=1}^mM_i\isom\oplus_{j=1}^nN_i$
are two decompositions with indecomposable summands, then $m=n$ and,
after suitable renumbering, $M_i\isom N_i$ for all $i$.
\end{thm}
\begin{proof}
\cite[Theorem 1.4.6]{Ben}.
\end{proof}

This theorem allows us to speak of ``the'' indecomposable direct summands of a given
module. To study modules in terms of these summands, we must introduce cores,
projective covers, and loop spaces.

\begin{defn}
The (isomorphism class of the) direct sum of the non-projective indecomposable summands of
a given representation $M$ is called \emph{the core of $M$}, and will be denoted by
$\core(M)$. If we have $M\isom\core(M)$, we call $M$ itself \emph{a core}. The (isomorphism
class of the) direct sum of the projective indecomposable summands is called the \emph{projective
part of $M$}.
\end{defn}

\begin{defn}\hfill\begin{enumerate}
\item A homomorphism of modules is called \emph{essential} if it is surjective and
its restriction to every proper submodule of its domain is not surjective.
\item A \emph{projective cover} of a module $M$ consists of a projective module
$P$ and an essential map $\pi: P\To M$.
\end{enumerate}
\end{defn}

\begin{thm}
Any module has a projective cover, which is again finitely generated and unique
up to (non-unique) isomorphism. The projective cover of a direct sum is the direct sum
of the individual projective covers.
\end{thm}
\begin{proof}
\cite[Chapter 14, Proposition 4]{Ser1}.
\end{proof}
 We may thus speak of ``the'' projective cover $P_G(M)$ of a module.

It is known that the number of isomorphism classes of simple modules is finite
\cite[Chapter 18, Corollary 3]{Ser1}. By contrast, there are in general infinitely many
isomorphism classes of indecomposable modules \cite[Theorem 4.4.4]{Ben}. However,
the \emph{projective} indecomposable modules are easily described by the following lemma.

\begin{thm}\label{projectiveindecomposables}
The operation ``projective cover'' induces a bijection between the set
$\Irr G$ of isomorphism classes of simple modules and the set of isomorphism
classes of projective indecomposable modules.
\end{thm}
\begin{proof}
\cite[Chapter 14, Corollary 1]{Ser1}.
\end{proof}

It follows from the above theorem that any module $M$ has a decomposition
$$M\isom\core(M)\oplus\bigoplus_{S\in\Irr G}P_G(S)^{\oplus b(M,S)}$$
for unique integers $b(M,S)\ge 0$. To know the isomorphism class
of $M$ is to know its core and to know the value of these integers.

The core of a module is the degree zero case of a concept of ``loop spaces''
developed to understand modules ``up to projectives''. Other authors
write $\Omega_G^0(M):=\core(M)$. We will need the degree one case:

\begin{defn}
Given a module $M$, its \emph{(first) loop space} is
$$\Omega_G(M):=\Omega_G^1(M):=\ker(P_G(M)\To M).$$
Recursively, we define $\Omega_G^i(M):=\Omega_G(\Omega_G^{i-1}(M))$ for $i>1$.
\end{defn}

What follows are some technical lemmas. The reader only interested in the statements
of our theorems now has the necessary notation, and may skip the rest of this subsection.

\begin{prop}\label{homcalc}
Given a module $M$ and a simple module $S$, we have $$\Hom_G(M,S)=\Hom_G(P_G(M),S).$$
\end{prop}
\begin{proof}
We apply the functor $\Hom_G(-,S)$ to the exact sequence
$$0\To\Omega_G(M)\arrover{i} P_G(M)\arrover{\pi_M} M\To 0$$
to get the exact sequence
$$0\To\Hom_G(M,S)\To\Hom_G(P_G(M),S)\arrover{i^*}\Hom_G(\Omega_G(M),S).$$
The lemma follows from the equation $i^*=0$. Assume that $i^*\neq 0$, then there exists a
nonzero map $f\in\Hom_G(\Omega_G(M),S)$ which factors through $P_G(M)$ as $f=Fi$,
for some $F\in\Hom_G(P_G(M),S)$. Since $S$ is irreducible, $f$ and $F$ must be surjective.
The map $$P_G(M)\arrover{(\pi_M,F)}M\oplus S$$ is still surjective. Thus
$\ker F\arrover{\pi_M} M$ is surjective. Since $F\neq 0$, i.e. $\ker F\subsneq P(M)$,
this is a contradiction to the fact that $\pi_M$ is essential.
\end{proof}

\begin{cor}
If we write \[P_G(M)\isom\bigoplus_{S\in\Irr G}P_G(S)^{b(P_G(M),S)},\]
then the multiplicities $b(P_G(M),S)$ are given by $b(P_G(M),S)=\dim_k\Hom_G(M,S)$.
\end{cor}
\begin{proof}
We fix $S\in\Irr G$ and calculate by means of the previous lemma:

\begin{equation*}\begin{split}
\Hom_G(M,S)&=\Hom_G(P_G(M),S)\\
&\isom\bigoplus_{T\in\Irr_G}\Hom_G(P_G(T),S)^{\oplus b(P_G(M),T)}\\
&=\bigoplus_{T\in\Irr_G}\Hom_G(T,S)^{\oplus b(P_G(M),T)}.
\end{split}\end{equation*}
By Schur's lemma, the dimension of $\Hom_G(T,S)$ is $0$ or $1$, depending
on whether $T$ and $S$ are isomorphic or not. Thus the corollary follows
by counting dimensions.
\end{proof}

\begin{prop}\label{projectivitycriteria}
Given a module $M$, the following are equivalent:\begin{enumerate}
\item $M$ is projective,
\item $M$ is injective,
\item $\core(M)=0$,
\item $\Omega_G(M)$ is projective, and
\item $\Omega_G(M)=0$.
\end{enumerate}
Furthermore, $\Omega_G(M)$ is always a core.
\end{prop}
\begin{proof}
The equivalence of (i) and (ii) follows from \cite[Prop. 1.6.2]{Ben}
and \cite[Prop. 3.1.2]{Ben}.
Clearly, (i) and (iii) are equivalent by definition.
Since $\Omega_G(M)=0$ if and only if $P_G(M)\To M$ is an isomorphism, (i) and (v)
are equivalent.
The equivalence of (iv) and (v) follows from the claim that $\Omega_G(M)$ is
a core, which we now prove:

Assume that $P\subset\Omega_G(M)\subset P_G(M)$ is a non-zero projective submodule.
Then (by the equivalence of (i) and (ii)) $P_G(M)$ decomposes as a direct
sum $P_G(M)\isom P\oplus Q$, and the image of $Q$ in $M$ is all $M$. This
is a contradiction to the fact that $P_G(M)\To M$ is essential; hence $\Omega_G(M)$
is a core.
\end{proof}

The following proposition is well-known; we give a proof here
since it will be a central component in the proof of our Theorem \ref{theoremA}.
\begin{prop}\label{corecalc}
Consider an exact sequence $0\To N\To P\To M\To 0$ of modules,
where $P$ is projective. Then there exists an isomorphism
\[\core(N)\isom\Omega_G(M).\]
Furthermore, if we denote the projective part of $N$ by $Q$,
we have $P\isom P_G(M)\oplus Q$.
\end{prop}
\begin{proof}
We construct the following commutative diagram
\[\begin{CD}
0 @>>> N @>>> P  @>>> M  @>>> 0\\
@. @VVV @VVV @| @.\\
0 @>>> \Omega_G(M) @>>> P_G(M)  @>>> M  @>>> 0\\
\end{CD}\]
The middle vertical arrow exists because $P$ is projective; it is surjective
because $P_G(M)\To M$ is
essential. Let $Q$ be the kernel of this middle arrow. Since $P_G(M)$ is
projective, so is $Q$. By the snake lemma, the first vertical arrow
is surjective, and its kernel is isomorphic to $Q$. Since $Q$ is injective
(Proposition \ref{projectivitycriteria}), we have an isomorphism $$N\isom\Omega_G(M)\oplus Q,$$
which proves the second claim.
Since $Q$ is projective, and $\Omega_G(M)$ is a core (Proposition \ref{projectivitycriteria}),
$\Omega_G(M)$ is the core of $N$.
\end{proof}

\begin{prop}\label{divisibility}
Let $p^n$ be the $p$-part of the order of $G$, i.e. $|G|=p^nk$ with $k\in\NN$ and $p\nmid k$.
Then the dimension of every projective module is divisble by $p^n$.
\end{prop}
\begin{proof}
\cite[Exercise 16.3]{Ser1}.
\end{proof}

\begin{prop}\label{invprojectives}
Let $N\subset G$ be a normal subgroup, and consider the group $H:=G/N$.
There is an inclusion $\Irr H\subset\Irr G$. Given $S\in\Irr G$, we have
\[P_G(S)^N\isom\begin{cases}
P_H(S) & \text{if $S\in\Irr H$}\\
0 & \text{if $S\in\Irr G\setminus\Irr H$}.
\end{cases}\]
\end{prop}
\begin{proof}
\cite[Lemma 2.7]{Bor}.
\end{proof}

In the last section we will need the following statement about
group cohomology.

\begin{prop}\label{psylow}
Let $K\subset G$ be a subgroup with $p\nmid[G:K]$. Then for any
representation $M$ of $G$, and for all $i\ge 0$, the restriction map
\[\Res: H^i(G,M)\To H^i(K,M)\]
is injective. In particular, if $p\nmid|G|$, then $H^i(G,M)=0$ for all $i>0$.
\end{prop}
\begin{proof}
\cite[Corollary 3.6.18]{Ben}.
\end{proof}

\section{The $p$-rank of curves}\label{section2.1}

We continue to assume as given an algebraically closed field $k$ of characteristic
$p>0$. In this article, a \emph{curve} signifies a complete, smooth, connected,
$1$-dimensional variety over $k$. The \emph{(absolute) Frobenius morphism $F$}
of such a curve $X$ is the (canonical)  morphism which is the identity on topological spaces,
and the $p$-power map on sections of the structure sheaf. It induces maps on the (Zariski) cohomology groups
$H^i(X,\OO_X)$. These are additive, but not $k$-linear maps: They are $p$-linear, meaning that
$$F(\lambda\xi)=\lambda^pF(\xi)\quad\text{for $\lambda\in k$ and $\xi\in H^i(X,\OO_X)$}.$$
The only nontrivial case for curves is the induced map on $H^1(X,\OO_X)$.

For this, let us review some material on $p$-linear maps. There is a category of $p$-linear maps,
with objects the pairs $(V,F)$ consisting of a finite-dimensional vector space $V$ and a $p$-linear
endomorphism $F$ of $V$. The morphisms in this category are the linear maps on the underlying vector spaces which
commute with the given $p$-linear endomorphisms. Given such an object $(V,F)$, we set
$V^F:=\set{v\in V: Fv=v}$, the fixed vectors of $F$ in $V$, furthermore $V^s:=\bigcap_{i>0}\im F^i$,
and $V^n:=\bigcup_{i>0}\ker F^i$.

The integer $h=\dim_k V^s$ is called the \emph{stable rank of $F$}. The vector space
$V^s$ is often called the \emph{semisimple} part of $V$.

\begin{prop}\label{plinearmaps}
In the above situation, we have\begin{enumerate}
 \item $V^F$ is a $\FF_p$-vector space.
 \item $V^s$ and $V^n$ are $k$-vector spaces stable under $F$.
 \item $\dim_kV^s=\dim_{\FF_p}V^F$.
 \item $V=V^s\oplus V^n$.
 \item $F$ restricted to $V^s$ is bijective, $F$ restricted to $V^n$ is nilpotent.
 \item $(-)^s$ is an exact functor on the category of $p$-linear maps.
\end{enumerate}
\end{prop}
\begin{proof}
See \cite{Has} or \cite{Ser3} for (i) to (v). The last statement is clear, since
we assume the maps in the category to be compatible with the respective $p$-linear
maps $F$.
\end{proof}

On the dual vector space $V^*=\Hom_k(V,k)$ we can define a map $C$
by setting $C(\psi)(v):=\psi(F(v))^p$ for $v\in V$ and $\psi\in V^*$. This map is additive and $1/p$-linear,
i.e. we have $C(\lambda\psi)=\lambda^{1/p}C(\psi)$. The decomposition $V=V^s\oplus V^n$
corresponds to a decomposition of $V^*$, and $C$ has the same stable rank as $F$.
Since any $1/p$-linear map can be viewed as the dual of a $p$-linear map, the structure
theory of the previous proposition can be translated to $1/p$-linear maps.

\begin{defn}
The $p$-rank $h_X$ of a curve $X$ is the stable rank
of the Frobenius morphism on $H^1(X,\OO_X)$.
\end{defn}

It is clear that we have estimates $0\le h_X\le g_X$, where $g_X=\dim_k H^1(X,\OO_X)$ is
the genus of $X$.

We would like to know explicitly the dual map of the Frobenius morphism on $H^1(X,\OO_X)$.
Recall that a rational function $t\in k(X)$ is called \emph{separating}
if the field extension $k(X)/k(t)$ is separable. Given a meromorphic
differential $\omega=f\cdot dt$, where $f\in k(X)$, we may write
$$f=f_0^p+f_1^pt+\cdots+f^p_{p-1}t^{p-1}.$$
The \emph{Cartier operator} on differentials is defined by setting
\begin{equation}\label{Cartierdef}
\CC(\omega):=f_{p-1}dt= \left(\sqrt[p]{-\left(\frac{d}{dt}\right)^{p-1}f}\right)dt.
\end{equation}
This is well-defined and independent of the choice of $t$ \cite{Ser3}.

\begin{prop}\label{dualit}
The dual vector space of $H^1(X,\OO_X)$ is, by Serre duality,
the vector space $H^0(X,\Omega_X)$. Under this identification,
the Cartier operator is the dual map of the Frobenius morphism.
\end{prop}
\begin{proof}
\cite{Ser3}.
\end{proof}

The geometric meaning of the $p$-rank is the following: There
are $p^{h_X}$ unramified Galois coverings of the curve $X$
with Galois group $\FF_p$, up to isomorphism
of the covering curve together with the action of $\FF_p$. More precisely,
the group $\Hom(\pi_1^{et}(X),\FF_p)$ classifies such
covers, and there are natural isomorphisms
$$(H^0(X,\Omega_X)^C)^*\isom H^1(X,\OO_X)^F\isom H^1(X_{et},\FF_p)
\isom\Hom(\pi_1^{et}(X),\FF_p),$$
compatible with the operation of automorphisms of $X$ on the respective
vector spaces. For proofs and further background, we refer to the survey
in \cite{Bou}. In this article, we will avoid rationality questions in
representation theory by studying $H^1(X,\OO_X)^s=H^1(X_{et},\FF_p)\otimes_{\FF_p} k$
instead of $H^1(X,\OO_X)^F=H^1(X_{et},\FF_p)$.
Also, we will study the dual representation $H^0(X,\Omega_X)^s$ instead of
$H^1(X,\OO_X)^s$ to simplify computations.



\section{The cores of $p$-rank representations}

Consider a curve $X$, and a finite subgroup $G\subset\Aut(X)$. If $X$ is of
genus $g_X\ge 2$, then $\Aut(X)$ itself is finite, but even in that case
we wish to allow ourselves the freedom of choosing a smaller group.

Given a point $x\in X$, we use the notation $v_x(-)$ for the function which
assigns to a function, differential or divisor its order at $x$.

\begin{prop}\label{commutingoperations}
Let $D$ be an effective divisor on a curve $X$. The Cartier operator $\CC$ operates
on the sheaf $\Omega_X(D)$. If $D$ is $G$-invariant, then $G$ also operates
on this sheaf, and the two operations commute.

In particular, the vector space $H^0(X,\Omega_X(D))^s$ of semisimple
differentials with respect to $\CC$ is a (finite dimensional)
representation of $G$.
\end{prop}
\begin{proof}
Consider an open set $U\subset X$ and a differential $\omega\in\Omega_X(D)(U)$.
For $x\in U$ choose a local parameter $t$ at $P$ and write
\begin{equation}\label{cschr}\omega=(f_0^p+f_1^p\cdot t+\cdots+f_{p-1}^p\cdot t^{p-1})dt = f\cdot dt\end{equation}
as in Section \ref{section2.1}, noting that $t$ is separating. Setting $n=v_x(D)\ge 0$, the assumptions
imply that $v_x(f)=v_x(\omega)\ge -n$. Thus the estimate
$$p\cdot v_x(f_{p-1})+p-1=v_x(f^p_{p-1}t^{p-1})\ge \min_i(v_x(f^p_it^i)=v_x(f)\ge -n$$
holds true. We now see that $v_x(\CC(\omega))=v_x(f_{p-1})\ge\lceil\frac{1-p-n}{p}\rceil\ge -n$,
where $\lceil y \rceil$ signifies the least entire number greater than $y$.
Therefore, we have $\CC(\omega)\in\Omega_X(D)(U)$.

Choose $g\in G$. We have $\CC(\omega)^g=(f_{p-1}dt)^g=f_{p-1}^g(dt)^g$.
On the other hand, if $t$ is separating, so is $s=t^g$, thus if we
write $\omega^g=(\cdots+(f_{p-1}^g)^p\cdot s^{p-1})ds$, we have
$\CC(\omega^g)=f_{p-1}^g ds=\CC(\omega)^g$, since the definition
of $\CC$ does not depend on the choice of separating variable.
\end{proof}

\begin{defn}
The module $V_D:=H^0(X,\Omega_X(D))^s$ of the previous proposition is the \emph{$p$-rank representation
of $G$ associated to the divisor $D$}.
\end{defn}

We introduce the notion $D^{red}$ for the reduced effective divisor
associated to $D$. The following observation will prove to be helpful:

\begin{prop}\label{DisDred}
If $D$ is an effective divisor on $X$, then the $p$-rank representation does
not depend on the multiplicities of $D$, i.e. $V_D=V_{D^{red}}.$
\end{prop}
\begin{proof}
The claim is that semisimple differentials have poles of order at most one.
By Proposition \ref{plinearmaps}(iii), it is sufficient to prove this claim for differentials
of the form $\omega=\CC(\omega)$. If $v_x(\omega)=-n<0$, then as in the
proof of Proposition \ref{commutingoperations}, we have $v_x(\CC(\omega))\ge\lceil\frac{1-p-n}{p}\rceil$.
It is elementary to prove that $$\frac{1-p-n}{p}=-n \iff n=1,$$
so we see that $v_x(\omega)\ge -1$.
\end{proof}

In the following, we will always assume that $D$ and $\Dtil$ are are effective
and $G$-invariant reduced divisors.

\begin{defn}
We will call $\Dtil$ \emph{sufficiently large with respect to $G$}
if it is non-empty, and contains all points of $X$ with nontrivial
stabilizer in $G$.
\end{defn}

\begin{prop}[Nakajima]\label{nakajimaprop}
If $\Dtil$ is sufficiently large with respect to $G$, the $p$-rank
representation $V_{\Dtil}$ is a projective $k[G]$-module.
\end{prop}
\begin{proof}
Let $P\subset G$ be a $p$-Sylow subgroup of $G$. By \cite[Theorem 1]{Nak2}
we know that $V_{\Dtil}$ is $k[P]$-free. This is equivalent
to the fact that $V_{\Dtil}$ is $k[G]$-projective \cite[Corollary 3.6.10]{Ben}.
\end{proof}

We will present the core of a $p$-rank representation
as a loop space of the following \emph{ramification module}.

\begin{defn}
Given a $G$-invariant effective reduced divisor $D$ as above, we \emph{choose} a sufficiently
large divisor $\Dtil\supset D$. The \emph{ramification module of $V_D$
(with respect to $\Dtil$}) is the following:

\[R_{G,D,\Dtil}:=\begin{cases}
k[\Dtil\setminus D], & \text{if $D\neq\emptyset$}\\
\ker(k[\Dtil]\To k,\quad \sum\lambda_xx\mapsto \sum\lambda_x), & \text{if $D=\emptyset$},
\end{cases}\]
where, for any reduced effective divisor $E$, by $k[E]:=\bigoplus_{x\in E}k\cdot x$ we denote
the affine coordinate ring of the reduced subvariety of $X$ associated to $E$.

The \emph{core module of $V_D$} is the loop space
$$C_D:=\Omega_G(R_{G,D,\Dtil}).$$
\end{defn}

\begin{rem}
We note that the module $C_D$ does not depend on the choice of
$\Dtil$, since enlarging $\Dtil$ corresponds to
adding to $R_{G,D,\Dtil}$ direct summands isomorphic to $k[G]$,
and such free summands are annihilated
by the loop space operator. Furthermore,
$$k[\Dtil\setminus D]\isom\bigoplus_{x\in\Dtil\setminus D\:(\MOD G)}k[G/G_x]$$
is a sum of induced representations of the trival representation.
\end{rem}

\begin{thm}\label{theoremA}
The core of the $p$-rank representation associated to a $G$-invariant
effective divisor $D$ (not necessarily reduced) is given by
the following formula:
$$\core(V_D)\isom C_{D^{red}}.$$
\end{thm}
\begin{proof}
By Proposition \ref{DisDred} we may assume that $D$ is reduced. We choose $\Dtil\supset D$
sufficiently large. Then $\Dtil\setminus D$ is also reduced and $G$-invariant,
and the residue map induces an exact sequence
$$0\To\Omega_X(D)\To\Omega_X(\Dtil)\arrover{\Res}\OO_{\Dtil\setminus D}\To 0,$$
which is invariant under the operation of $G$. It induces a long exact sequence
\begin{eqnarray*}0&\To& H^0(X,\Omega_X(D))\To
H^0(X,\Omega_X(\Dtil))\To k[\Dtil\setminus D]\\
&\arrover{\delta} &H^1(X,\Omega_X(D))\To 0,
\end{eqnarray*}
which terminates at $H^1(X,\Omega_X(\Dtil))=0$ since $\Dtil\neq 0$. Clearly, we have
$\ker\delta = R_{G,D,\Dtil}$; hence there is an exact sequence of $k[G]$-modules
$$0\To H^0(X,\Omega_X(D))\To H^0(X,\Omega_X(\Dtil))\arrover{\Res} R_{G,D,\Dtil}\To 0.$$
In order to extract from this an exact sequence of semisimple parts, we define a
$1/p$-linear map on $k[\Dtil\setminus D]=\bigoplus_{d\in\Dtil\setminus D}k\cdot d$
by letting it operate on the standard basis $\set{d}_{d\in\Dtil\setminus D}$ as the identity.
This induces a $1/p$-linear map on $R_{G,D,\Dtil}$, compatible
with the operation of $G$. Since we know that $\Res(\CC\omega)^p=\Res(\omega)$ by \cite{Ser3},
the above sequence is an exact sequence in the category of $1/p$-linear maps.
Thus, by the exactness of $(-)^s$, we obtain the exact sequence
$$0\To V_D\To V_{\Dtil}\To R_{G,D,\Dtil}\To 0.$$
By Proposition \ref{nakajimaprop} the middle term is a projective
module, and Proposition \ref{corecalc} gives the desired result.
\end{proof}

\begin{rem}
If $G$ has no fixed points and $D=\emptyset$,
the core of the associated $p$-rank representation is
$$\core(V_\emptyset)=\Omega^1_G(R_{G,\emptyset,\Dtil})=\Omega^2_G(k),$$
since $k[\Dtil]\isom k[G]^r$ for some $r\ge 1$, which implies that
$\core(R_{G,\emptyset,\Dtil})=\Omega^2_G(k)$. This particular core
has been calculated by Borne in \cite{Bor}.
\end{rem}

\begin{rem}
Since a projective representation is determined up to isomorphism by its
composition factors  \cite[Chapter 14, Corollary 3]{Ser1}, the
local invariants used in Theorem \ref{theoremA} and the modular character of a $p$-rank
representation determine such a representation up to \emph{isomorphism}.
\end{rem}

\begin{cor}
Consider a curve $X$ and a finite group $G$ of automorphisms of $X$.
Let $r$ be the number of fixed points of $G$ on $X$,
and let $p^n$ be the $p$-part of the order of $G$. Then:
$$h_X\equiv 1-r \pmod{p^n} $$
\end{cor}
\begin{proof}
We choose a minimal sufficiently large divisor $\Dtil\supset\emptyset$,
and set $R:=R_{G,\emptyset,\Dtil}$.
Since $h_X$ is the dimension of $V_{\emptyset}$
and by Theorem \ref{theoremA} this module
differs from its core only by projective summands, Proposition \ref{divisibility}
implies that \[h_X\equiv\dim\Omega_G(R) \pmod{p^n}.\]
Similar reasoning applies to $\Omega_G(R)$ and $R$,
which have dimensions adding
up to the dimension of the projective module $P_G(R)$, and shows that
\[\dim\Omega_G(R)\equiv-\dim R \pmod{p^n}.\]
If $G$ has fixed points,
then $\dim R=\dim k[\Dtil]-1=r-1$, if $G$ has no fixed points, then
$\dim R=\dim k[\Dtil]-1=p^n-1\equiv r-1\pmod{p^n}$;
hence we can combine the above congruences to obtain the corollary.
\end{proof}


\section{Borne invariants of quotient curves}

In addition to the notation and conventions of the previous section,
we consider a normal subgroup $N$ of $G$, and the short exact sequence
$$1\To N\To G\To H\To 1.$$A representation of $H$
lifts to a representation of $G$, and we obtain an inclusion $\Irr H\subset\Irr G$
of the set of irreducible representations.

Let $Y:=X/N$ be the quotient curve, and let $\pi: X\To Y$ be the canonical
projection. There is a natural induced operation of $H$ on $Y$. The notation
of the last section will sometimes have to be decorated by subscripts $X$ or $Y$.

\begin{defn}
The \emph{Borne invariants $b(G,D,S)$ of the curve $X$ (with respect to $G$ and $D$)}
are the multiplicities of the projective indecomposable modules in the $p$-rank
representation of $G$ with respect to $D$. Thus, we have an isomorphism
$$V_D=\core(V_D)\oplus\bigoplus_{S\in\Irr G}P_G(S)^{\oplus b(G,D,S)}.$$
We simplify notation, setting $b(G,S):=b(G,\emptyset,S)$.
\end{defn}

\begin{prop}[Pink]\label{calcinv}
Let $D$ be an $N$-invariant reduced effective divisor on $X$.
There is a natural isomorphism of sheaves
\[\pi_*\Omega_X(D)^N\isom\Omega_Y(E)\] for an
effective divisor $E$ on $Y$, which commutes with the Cartier operator and the operation
of $G$. We have
\begin{equation*}
  E^{red}=\pi(D)^{red}\cup \set{y\in Y\mid \pi\text{ is wildly ramified over $y$}}.
\end{equation*}
\end{prop}
\begin{proof}
Pulling back differentials from $Y$ to $X$ via $\pi$ induces an injective sheaf homomorphism
\[\Omega_Y\To\pi_*\Omega_X(D)^N.\]
The target of this homomorphism is a torsion-free, coherent sheaf of rank $1$; hence there
is a unique effective divisor $E$ on $Y$ such that the above homomorphism extends to
an isomorphism $\Omega_Y(E)\To\pi_*\Omega_X(D)^N$.
By construction of the Cartier operator, this is Cartier-equivariant.

We now proceed to determine $E$. If $R$ is a local ring, we denote its completion by $\widehat{R}$.
Choose $y\in Y$, we then have
\begin{equation*}\begin{split}
\pi_*\Omega_X(D)^N\otimes_{\OO_{Y,y}}\widehat{\OO_{Y,y}}=&\left(\bigoplus_{x\in\pi^{-1}(y)}\Omega_X(D)\otimes_{\OO_{X,x}}\widehat{\OO_{X,x}}\right)^N\\
=&\left(\Omega_X(D)\otimes_{\OO_{X,x}}\widehat{\OO_{X,x}}\right)^{N_x}\quad\text{for any $x\in\pi^{-1}(y)$.}\\
=&\:\widehat{\Omega_X(D)_x}^{N_x}
\end{split}\end{equation*}
Choose $x\in\pi^{-1}(y)$, and denote again by $x$ and $y$ local parameters at $x$
and $y$ respectively. We have $\widehat{\OO_{X,x}}=k[[x]]$ and
$\widehat{\OO_{Y,y}}=k[[y]]=k[[x]]^{N_x}$. Setting $n:=|N_x|$, we may express
$y$ as $$y=x^n + \text{terms of higher order in $x$}.$$
It follows that $m:=v_x(\frac{dy}{dx})=n-1$ if $p\nmid n$, and $m\ge n$ if $p\mid n$.
Let us set $d:=v_x(D)$ and write $d+m=an+b$, for integers $a,b\ge 0$ with $b\le n-1$.
The following are equivalent:
\begin{enumerate}
  \item $a\ge 1$
  \item $d+m\ge n$
  \item $d\ge 1$ or $p \mid n$
\end{enumerate}

We now see that
\[\widehat{\Omega_X(D)_x}^{N_x}=\left(\frac{1}{x^d}k[[x]]dx\right)^{N_x}=(\frac{1}{x^{d+m}}k[[x]])^{N_x}dy=\frac{1}{y^a}k[[y]]dy,\]
which implies that $\pi_*\Omega_X(D)^N=\Omega_Y(E)$ if we set $v_Q(E):=a$, as claimed.
\end{proof}

\begin{defn}
Consider an irreducible representation $S\in\Irr G$. The restriction maps $H^1(G,S)\To H^1(G_x,S)$
combine to a global restriction map\begin{equation*}
  r_{G,X,S}: H^1(G,S)\To\prod_{x\in X}H^1(G_x,S)=\bigoplus_{x\in X^G}H^1(G_x,S).
\end{equation*}
Its kernel $H^1_{LT,X}(G,S):=\ker r_{G,X,S}$ consists
of the \emph{locally trivial first cohomology classes of S} (with respect
to $G$ and $X$). We set \[d(G,X,S):=\dim H^1_{LT,X}(G,S).\]
\end{defn}

\begin{thm}\label{theoremB}
The Borne invariants of $X$ and $Y=X/N$ with respect to $G$ and $H=G/N$
for $T\in\Irr H$ are related by the following formula:\begin{equation*}
b(G,T)+d(G,X,T)=b(H,T)+d(H,Y,T).
\end{equation*}
\end{thm}
\begin{proof}
This is a lengthy calculation, which we divide into several steps.
We choose a sufficiently large divisor $\Dtil$ on $X$ with respect to $G$, and
set $\Etil:=\pi_*(\Dtil)^{red}$; this is a sufficiently large divisor on $Y$ with respect
to $H$.

\emph{Step 1:} Since $\Dtil$ contains all ramified points, wild or not,
Proposition \ref{calcinv} implies that
\begin{eqnarray*}
(V_{X,\Dtil})^N &=& V_{Y,\Etil}.
\end{eqnarray*}
In particular, since $V_{X,\Dtil}$ is projective by Proposition \ref{nakajimaprop},
we may apply Proposition \ref{invprojectives} to its indecomposable summands
to obtain \begin{equation}\label{e1}
b(G,\Dtil,T)=b(H,\Etil,T)\quad\text{for $T\in\Irr H$}.
\end{equation}

\emph{Step 2:} The short exact sequence
$$0\To V_{X,\emptyset}\To V_{X,\Dtil}\To R_{G,\emptyset,\Dtil}\To 0$$
induces, by the second claim of Proposition \ref{corecalc}, an isomorphism
$$\bigoplus_{S\in\Irr G}P_G(S)^{b(G,\Dtil,S)}\isom P_G(R_{G,\emptyset,\Dtil})\oplus\bigoplus_{S\in\Irr G}P_G(S)^{b(G,\emptyset,S)}.$$
In particular, using Proposition \ref{homcalc}, we may apply $\Hom_G(-,S)$ to
deduce the equation \begin{equation}\label{e2}
b(G,\Dtil,S)=\dim_k\Hom_G(R_{G,\emptyset,\Dtil},S)+b(G,\emptyset,S)\quad\text{for $S\in\Irr G$}.
\end{equation}

\emph{Step 3:} On the other hand, let us consider $S\in\Irr G$ and the short exact sequence
$$0\To R_{G,\emptyset,\Dtil}\To k[\Dtil]\To k\To 0.$$
Applying $\Hom_G(-,S)$ to this sequence gives an exact sequence\begin{eqnarray*}
0&\To &S^G\To\bigoplus_{x\in\Dtil\:(\MOD G)}S^{G_x}\To\Hom_G(R_{G,\emptyset,\Dtil},S)\\
 &\To &H^1(G,S)\arrover{r_{G,X,S}}\bigoplus_{x\in\Dtil\:(\MOD G)}H^1(G_x,S),
\end{eqnarray*}
that is, an exact sequence\[
0\To S^G\To\bigoplus_{x\in\Dtil\:(\MOD G)}S^{G_x}\To\Hom_G(R_{G,\emptyset,\Dtil},S)
 \To H^1_{LT,X}(G,S)\To 0.\]
Similar reasoning applies to $T\in\Irr H$, leading to the exact sequence\[
0\To T^H\To\bigoplus_{y\in\Etil\:(\MOD H)} T^{H_y}\To\Hom_H(R_{H,\emptyset,\Etil},T)
 \To  H^1_{LT,Y}(H,T)\To 0.\]
Using the fact that the alternating sum of dimensions in an exact
sequence is $0$, the equality $T^G=T^H$ and, for $y=\pi(x)$, the similar
equalities $V^{G_x}=V^{H_y}$,
we have \begin{equation}\label{e3}
\dim_k\Hom_H(R_{H,\emptyset,\Etil},T)-\dim_k\Hom_G(R_{G,\emptyset,\Dtil},T)=d(H,Y,T)-d(G,X,T).
\end{equation}

\emph{Step 4:} Finally, combining equations \eqref{e1} and \eqref{e2} (for $X$ \emph{and} $Y$),
and \eqref{e3} gives the result.
\end{proof}

\begin{rem}
If $N$ is a $p$-group, then it is known that $\Irr G=\Irr H$ \cite[Remark after Definition 2.5]{Bor}.
Thus, in this case, the Borne invariants of $Y$ determine \emph{all} the Borne invariants
of $X$. In this sense, Theorem \ref{theoremB} generalizes the equivariant Deuring-Shafarevich formula
of Shoichi Nakajima \cite{Nak2}, which is the special case of $N=G$ being a $p$-group.
\end{rem}

\begin{rem}
If the operation of $G$ on $X$ is \emph{tame}, that is if $p\nmid|G_x|$ for all $x\in X$,
then all higher cohomology groups of the stabilizers $G_x$ vanish by Proposition \ref{psylow}.
Thus $d(G,X,S)=\dim H^1(G,S)$, which proves the
conjecture that Niels Borne states in \cite{Bor} after Proposition 2.4.
\end{rem}

Under certain circumstances, the calculation of the locally trivial cohomology
groups is not necessary:

\begin{prop}
The following estimate holds true: $$b(G,T)\le b(H,T)\quad\text{for all $T\in\Irr H$}.$$
Furthermore, if there is an $x\in X$ such that $p\nmid[G:G_x]$ or $p\nmid[N:N_x]$, then
\[b(G,T)=b(H,T)\quad\text{for all $T\in\Irr H$}.\]
\end{prop}
\begin{proof}
Given $x\in X$ and setting $y=\pi(x)$, the sequence $$1\To N_x\To G_x\To H_y\To 1$$is
exact. We choose $T\in\Irr H$ and use the abbreviations $LT_G:=H^1_{LT,X}(G,T)$
and $LT_H:=H^1_{LT,Y}(H,T)$. The inflation-restriction sequence of group
cohomology \cite[Chapter 3.4, Exercise]{Ben} gives the exact sequence
$$0\To H^1(H,T)\arrover{\inf} H^1(G,T)\arrover{\res} H^1(N,T)^H.$$We use this to construct the
following commutative diagram:

\[\begin{CD}
    @.    0     @.     0        @.            0\\
@.       @VVV          @VVV                   @VVV\\
0 @>>> LT_H     @>>>  H^1(H,T)   @>a>> \bigoplus_y H^1(H_y,T)\\
@.      @ViVV          @VVV                   @VVV\\
0 @>>> LT_G     @>>>  H^1(G,T)   @>b>> \bigoplus_x H^1(G_x,T)\\
@.       @VVV          @VVV                   @VVV\\
  @.   \coker i @>j>> H^1(N,T)^H @>c>> \bigoplus_x H^1(N_x,T)^{H_y}\\
@.      @VVV             @.                     @.\\
@. 0 @. @. @.
\end{CD}\]

The first two rows are exact by the definition of locally trivial cohomology classes.
The last two columns are exact by the inflation-restriction sequence.
The injectivity of the inflation maps shows that $i$ is injective, so the
first column is exact and in particular $d(H,Y,T)\le d(G,X,T)$,
which implies that $b(G,T)\le b(H,T)$.

Now a diagram chase shows that $LT_H=LT_G\cap H^1(H,T)$; hence the
induced map $j$ is injective and its image lies in $\ker c$. To show
that $b(G,T)=b(H,T)$ in the cases mentioned in the proposition, we will
show that $LT_G=LT_H$.

If $p\nmid[G:G_x]$ holds for some $x\in X$, then by Proposition \ref{psylow} the restriction
map $H^1(G,T)\To H^1(G_x,T)$ is injective. Thus, $b$ is injective,
and $LT_G=\ker b = 0$. Since $LT_H\subset LT_G=0$,
it follows that $LT_H=LT_G=0$.

If $p\nmid[N:N_x]$ holds for some $x\in X$, then by Proposition \ref{psylow}
the map $H^1(N,T)\To H^1(N_x,T)$
is injective, so $c$ is injective. Since $\coker i\subset\ker c = 0$,
it follows that $LT_G=LT_H$.

\end{proof}


\bibliographystyle{amsalpha}

\end{document}